\documentclass[a4paper,12pt]{article}
\usepackage{amssymb}
\usepackage{amsmath}
\usepackage{graphicx}                 

\begin{document}

\title{Algorithmic aspects of algebraic methods\\
for graph isomorphism testing}

\author{M. Kauki\v c, Dept. of Mathematical Methods, FRI \v ZU, \\
Ve\v lk\'y Diel, 01026 \v Zilina, Slovakia\\
e-mail: {\it mike@frcatel.fri.utc.sk}}

\newtheorem{theorem}{Theorem}[section]
\newtheorem{corollary}[theorem]{Corollary}
\newtheorem{lemma}[theorem]{Lemma}
\newtheorem{exmple}[theorem]{Example}
\newtheorem{defn}[theorem]{Definition}
\newtheorem{proposition}[theorem]{Proposition}
\newtheorem{conjecture}[theorem]{Conjecture}
\newtheorem{rmrk}[theorem]{Remark}

\maketitle


\begin{abstract}
We present the implementation of an algorithm for graph isomorphism testing, based on ideas about number of walks (of sufficiently large length) between vertices, expanded for strongly regular graphs (SRG-s) by testing the local complements 
and values of determinants of their adjacency matrices. All known non-isomorphic SRG-s (with no more than 64 vertices) are distinguishable by this method. 
\end{abstract}

\noindent
{\bf Keywords}: {\it graph isomorphism, adjacency matrix, strongly regular graphs, 
local complement, Python, SciPy}

\section{Introduction, basic algorithm}

The graph isomorphism problem is one of difficult central problems in graph theory. It can be easily formulated as the question, whether two differently looking images represent the same graph. We will consider only undirected graphs without loops and multiple edges. Let us denote by $V(G), E(G)$ the sets of vertices and edges of graph $G$. Two graphs $G, H$ are isomorphic iff there exists a bijective mapping $f: V(G)\mapsto V(H)$ between sets of vertices of graphs $G$ and $H$ such that for all pairs of vertices $u,v \in V(G)$ the edge
$\{u,v\}$ belongs to $E(G)$ if and only if $\{f(u),f(v)\}$ is an edge
of graph $H$.

The {\it graph invariant} $\lambda$ is a mapping from the class of all graphs 
${\mathcal G}$ to the set ${\mathcal R}$ such that for any isomorphic graphs $G,H$
the images $\lambda(G), \lambda(H)$ are the same, i.e. 
$\lambda(G) = \lambda(H)$. Commonly used graph invariants are e.g. number of
vertices, number of edges, degree-sequence (the non-descending sequence of
degrees of the vertices), eigenvalues or characteristic polynomial of the adjacency matrix, etc.
All above-mentioned invariants are {\it incomplete}, which means that we can
find pairs of nonisomorphic graphs with the same invariants. 

There are also {complete graph invariants}, i.e. $\lambda(G) = \lambda(H)$ 
holds if and only if graphs $G, H$ are isomorphic (see  e.g. \cite{chstring}
for greatest characteristic string invariant). Such invariants are not known to be computable in polynomial time, because the existence of polynomial algorithm for the class of graph isomorphism problems is still an open question. The existing polynomial algorithms are inexact in the sense that they cannot recognize some pairs of nonisomorphic graphs (because they use incomplete invariants for testing). Our proposed algorithm also belongs to this type of non-exact polynomial algorithms.

Let us briefly summarize the main results from Czimmermann \cite{cimo2},
\cite{cimo1} where the basic variant of isomorphism testing algorithm using walks between vertices counting was formulated. 

Let $A=(a_{ij})$ be the adjacency matrix of graph $G$. Then the number of all possible walks of length $k$ from vertex $v_i\in V(G)$ to vertex $v_j\in V(G)$
is simply the element $a_{ij}^{(k)}$ of matrix power $A^k$. 
Let us assign to every vertex $v_i \in V(G)$  the set 
\begin{equation}
S_i=\{s_{i1},s_{i2},\dots,s_{in}\},\quad s_{ij}=\left\{{a_{ij}^{(k)}}\right\}_{k=1}^{m},
\label{nwalks}
\end{equation}
i.e. $s_{ij}$ is the sequence formed from number of walks 
of length $1,2,\dots, m$ which begin in vertex $v_i$ and end in $v_j$. 
Let us have a pair of graphs $G, H$ with equal number $n$ of vertices for
which the sequences of sets (\ref{nwalks}) are    
$${\mathcal W}_G=\{S_1,S_2,\dots S_n\},\  {\mathcal W}_H=\{T_1,T_2,\dots,T_n\}.$$ 
If $G, H$ are isomorphic, then for suitable permutation
of vertices of graph $G$ the two sequences are equal. The proposed algorithm 
is based on comparison of permuted sequences ${\mathcal W}_G, {\mathcal W}_H$ with
$m$ (maximal length of walks considered) sufficiently large. In \cite{cimo1}
it was shown that $m$ can be chosen as the number of common distinct
eigenvalues of graphs $G,H$ (there is no additional information coming
from walks of greater length).  
\begin{figure*}[h]
\begin{center}
\includegraphics[width=4.5cm]{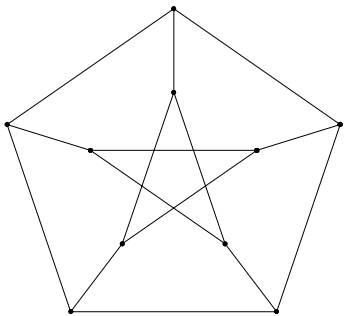}
\end{center}
\caption{Petersen graph -- strongly regular graph with parameters (10,3,0,1)}
\label{petersen}
\end{figure*}
For this basic form of algorithm it is easy to find the causes when it fails; e.g. two nonisomorpic {\em strongly regular graphs} with the same set of parameters are indistinguishable. The graph $G$ is called {\it strongly
regular (SRG) with parameters} $n,d, \alpha, \beta$  if it has $n$ vertices,
every vertex has degree $d$, every pair of adjacent vertices has $\alpha$
common neighbours and every pair of distinct nonadjacent vertices
has $\beta$ common neighbours. For example, the well-known Petersen graph
(see Fig. \ref{petersen}) is SRG with parameters (10,3,0,1). 

Each SRG has only three distinct eigenvalues, which are completely determined
by parameters $n,d, \alpha, \beta$. Thus, the complete information about
walks is contained in sequences of sets of walks with length $k \le 3$.
Such short walks are not sufficient to discover non-isomorfism of SRG-s
with the same set of parameters. For SRG-s, we need to investigate some additional properties.  

\section{Improved algorithm \\
for strongly regular graphs}

The key idea is to destroy the regularity of given SRG by representing it
with the sequence of suitable modified non-regular graphs. We will explore
the notion of {\it local complement} (see \cite{godsil}) for this purpose.

\begin{defn}
Let $G$ be a graph. The local complement $G_L(u)$ of $G$ at the vertex 
$u \in V(G)$ is defined to be the graph with the same vertex set as $G$ such that:
\begin{enumerate} 
\item 
if $v$ and $w$ are distinct neighbours of $u$, then they are adjacent
in  $G_L(u)$ if and only if they are not adjacent in G, 
\item 
if $v$ and $w$ are distinct vertices of $G$, and not both neighbours of $u$, then they are adjacent in $G_L(u)$ if and only if they are adjacent 
in $G$. 
\end{enumerate}
\end{defn}

Now, instead of testing isomorphism of two (SRG) graphs $G,H$, we can test for
existence of suitable permutation of vertices of graph $H$ such that the sequences of local complements of $G,H$ have the same graph invariants (e.g. the sequences of sets of walks, mentioned in previous section). 
%
%
\begin{table}[h]
$$
\begin{minipage}{5cm}
\begin{tabular}{|c|c|c|c|c|}
\hline
n & d & $\alpha$ & $\beta$ & No.\\
\hline\hline
16 & 6   & 2  & 2  & 2 \\
\hline
25 & 12  & 5  & 6  & 15 \\
\hline
26 & 10  & 3  & 4  & 10 \\
\hline
28 & 12 & 6 & 4 & 4 \\
\hline
29 & 14 & 6 & 7& 41 \\
\hline
35 & 16 & 6 & 8 & 3854 \\
\hline
35 & 18 & 9 & 9 & 227 \\
\hline
\end{tabular}
\end{minipage}
\begin{minipage}{5cm}
\begin{tabular}{|c|c|c|c|c|}
\hline
n & d & $\alpha$ & $\beta$ & No.\\
\hline\hline
36 & 14 & 4 & 6 & 180 \\
\hline
36 & 15 & 6 & 6 & 32548 \\
\hline
37 & 18 & 8 & 9 & 6760 \\
\hline
40 & 12 & 2 & 4 & 28 \\
\hline
45 & 12 & 3 & 3 & 78 \\
\hline
64 & 18 & 2 & 6 & 167 \\
\hline
\end{tabular}
\vspace{\baselineskip}
\end{minipage}
$$
\caption{Groups of nonisomorfic SRG-s}
\label{srgtab}
\end{table}

In Table \ref{srgtab} there are parameters and numbers of nonisomorphic
SRG-s with given parameter sets for up to $64$ vertices. We obtained
adjacency matrices for all mentioned graphs from WEB-sources \cite{srgdb1}, \cite{srgdb2}. To make the later processing easy, the SQL database of SRG-s
were created with the aid of simple but powerful database engine
{\it SQLite} (see \cite{sqlite}) and the corresponding interface module 
{\it pysqlite} \cite{pysqlite} for programming language Python. 

In the paper \cite{qwalks} the authors were able to distinguish all
graphs from this table (although they have omitted two large groups
of SRG-s with parameters $(35,16,6,8)$ and $(37,18,8,9)$) with their
algorithm based on eigenvalues of certain matrix inspired by the notion of 
coined quantum walks. We will show that our modified algorithm can
distinguish all of the graphs in Table \ref{srgtab}, when combined
with such easy to compute incomplete graph invariants as determinants
of adjacency matrices of local complements.  

\section{Algorithm implementation and testing}

We will focus on testing our algorithms only for set of SRG-s 
from Table \ref{srgtab}, because strongly regular graphs cause 
main difficulties for graph isomorphism algorithms.

For the rapid implementation of basic algorithm
(see section 1) we choose the interpreted programming language {\it Python}
\cite{python}. It has very clear syntax and together with {\it SciPy}
(Scientific tools for Python,\cite{scipy})  and IPython (enhanced Interactive Python shell, see \cite{ipython}) it forms very powerful and convenient
environment for rapid development of prototype applications. 

In the basic algorithm we need to compute the powers $A^k$ of 0-1 matrices
$A$ (for SRG-s, up to $k=3$). Underlying C and FORTRAN libraries in SciPy
use machine integer and floating point numbers with fixed precision.
This will give wrong results when the elements of $A^k$ will be large
(e.g. the $17$-th power of adjacency matrix of the unique SRG 
with parameters $(36, 10, 4, 2)$ shows on 32-bit computers some negative entries, which is completely wrong). The same problem exists with computing determinants (for the above matrix even on 64-bit computer
we see the nonsensical value -351843720888319.81 of determinant). 

Python supports arbitrary precision integers, but it has no efficient linear
algebra operations with integer matrices and vectors. So we decided to 
complement the software tools used with NTL \cite{ntl}, a C++ library for manipulating arbitrary length integers, and for vectors, matrices, and polynomials over such integers. We used the {\it ctypes} Python module
(see \cite{ctypes}) to call the functions in  NTL library in Python.

We verified all ``small groups'' (with at most 227 members) from Table \ref{srgtab}, using only basic algorithm applied to local complements.
But the testing for three ``big'' groups of SRG-s was very time-consuming,
so we considered to make preliminary classification of nonisomorphic
graphs with the same parameters, using some cheap graph invariants.

In Czimmermann \cite{cimo3} there is shown that if two graphs $G,H$
have distinct determinants of their corresponding adjacency matrices, 
then the spectra of adjacency matrices are also distinct and
that two graphs with distinct spectra can be distinguished by the basic
variant of algorithm, using the sets of walks between vertices. 

The algorithm used for testing of three big groups
of SRG-s, i.e. the groups with parameters $(35, 16, 6, 8)$, 
$(36, 15, 6, 6)$ and $(37, 18, 8, 9)$ can be informally described
as follows:
\begin{enumerate}
\item {\it generate the sets of local complements 
${\mathcal S}_G = \{G_L(u), u \in V(G)\}$ for all graphs $G$ in the
group of nonisomorphic graphs with the same parameters}

\item {\it compute the determinants of adjacency matrices of graphs in
${\mathcal S}_G$, order the sequence of determinants by magnitude; 
next we partition the group of nonisomorphic SRG-s 
into the equivalence classes of graphs with the same 
(ordered) sequence of determinants of local complements; the graphs
belonging to a class containing more than one member are indistinguishable
by this graph invariant} 

\item {\it now it is sufficient to apply the basic variant of algorithm
only for testing isomorphism of pairs of graphs, contained inside the
same equivalence class.
}
\end{enumerate}

The group of 3854 graphs with parameters $(35, 16, 6, 8)$ 
has $44$ equivalence classes with more than one member, 
$42$ of them have exactly two members and there are two classes with four members. In the group of  32548 graphs with parameters $(36, 15, 6, 6)$ 
we found $160$ equivalence classes ($152$ with two members, $6$ with three
members and $2$ with four members). For the last group of 6760 graphs with
parameters $(37, 18, 8, 9)$ there are $3379$ equivalence classes, each with
two members. The basic algorithm applied to such small groups of graphs
is really computationally inexpensive.

\section{Conclusions}
The presented algorithm for isomorphism testing can efficiently recognize
all known small nonisomorphis SRG-s, thus there is a hope it will perform
well also on bigger strongly regular graphs. This needs further investigation.



\end{document}